\documentclass{amsart}

\newif\ifdraft\draftfalse
\newif\ifcite\citefalse
\newif\ifblow\blowtrue

\usepackage{amsfonts,amssymb, amsmath}
\usepackage{oldgerm}
\usepackage{amscd}
\usepackage{eucal} 
\usepackage{mathrsfs} 
\usepackage[hyperindex]{hyperref}  
\usepackage[alphabetic]{amsrefs}
\ifdraft\ifcite\usepackage{showkeys}\else\usepackage[notcite,notref]{showkeys}\fi\fi

\usepackage{xypic}
\newdir{ >}{{}*!/-5pt/\dir{>}}

\newtheorem{theorem}[equation]{Theorem}
\newtheorem{claim}[equation]{Claim}
\newtheorem{lemma}[equation]{Lemma}

\theoremstyle{remark}

\theoremstyle{definition}
\newtheorem{definition}[equation]{Definition}

\theoremstyle{remark}

\newtheorem{remark}[equation]{Remark}

\numberwithin{equation}{section}

\def\bp{\mathbb P}

\def\bean{\begin{eqnarray}}
\def\eean{\end{eqnarray}}
\def\bea{\begin{eqnarray*}}
\def\eea{\end{eqnarray*}}

\def\x0{X_{s_0}}
\def\2p{\bp^1\times \bp^1}

\def\a{\alpha}

\def\nm{\nonumber}
\def\er{\eqref}

\def\vol{d\sigma}

\def\ben{\begin{equation}}
\def\een{\end{equation}}

\newcommand\wt{\widetilde}

\newcommand\ind{\mathop{\mathrm{Ind}}\nolimits}
\newcommand\ppv{\partial_+ V}
\newcommand\pmv{\partial_- V}
\newcommand\vn{\vec n}
\newcommand\sing{\mathop{\mathrm{Sing}}\nolimits}
\newcommand\m{[\frac{n-1}{2}]}

\begin{document}
\title{On Sha's secondary Chern-Euler class}
\author{Zhaohu Nie}
\email{znie@psu.edu}
\address{Department of Mathematics\\
Penn State Altoona\\
3000 Ivyside Park\\
Altoona, PA 16601, USA}

\date{\today}

\begin{abstract}
For a manifold with boundary, the restriction of Chern's transgression form of the Euler curvature form over the boundary is closed. Its cohomology class is called the secondary Chern-Euler class and used by Sha to formulate a relative Poincar\'e-Hopf theorem, under the condition that the metric on the manifold is locally product near the boundary. 
We show that the secondary Chern-Euler form is exact away from the outward and inward unit normal vectors of the boundary by explicitly constructing a transgression form. Using Stokes' theorem, this evaluates the boundary term in Sha's relative Poincar\'e-Hopf theorem in terms of more classical indices of the tangential projection of a vector field. This evaluation in particular shows that Sha's relative Poincar\'e-Hopf theorem is equivalent to the more classical Law of Vector Fields. 
\end{abstract}

\maketitle

\section{Introduction}

Let $X$ be a smooth oriented compact Riemannian manifold with boundary $M$. Throughout the paper we fix $\dim X=n\geq 2$ and hence $\dim M=n-1$. 
On $M$, one has a canonical decomposition
\begin{equation}
\label{collar}
TX|_M\cong \nu\oplus TM,
\end{equation}
where $\nu$ is the rank 1 trivial normal bundle of $M$. 

In his famous proof \cite{chern,chern2} of the Gauss-Bonnet theorem, Chern constructed a differential form $\Phi$ (see \er{phi}) of degree $n-1$ 
on the tangent sphere bundle $STX$, consisting of unit vectors in $TX$, satisfying the following two conditions:
\begin{equation}
\label{dP=O}
d\Phi=-\Omega,
\end{equation}
where $\Omega$ is the Euler curvature form of $X$ (pulled back to $STX$) when $\dim X$ is even and $0$ 
otherwise, 
and 
\begin{equation*}
\wt{\Phi_0}=\wt{\vol}_{n-1},
\end{equation*}
i.e., the 0th term $\wt{\Phi_0}$ of $\Phi$ is the relative unit volume form for the fibration $S^{n-1}\to STX \to X$
(see \er{0th term}). 

By \er{dP=O}, one has
\ben
\label{d=0}
d\Phi=0\text{ on }STX|_M,
\een
since even if $\dim  X$ is even, $\Omega|_M=0$ by dimensional reason.
Following \cite{sha}, $\Phi$ on $STX|_M$ is called the \emph{secondary Chern-Euler form}, whose cohomology class is called the secondary Chern-Euler class. 

Secondary Chern-Euler classes are useful in studying the relative Poincar\'e-Hopf theorem. Let $V$ be a smooth vector field on $X$.
We assume that $V$ has only isolated singularities, i.e., the set $\sing V:=\{x\in X|V(x)=0\}$ is finite, and that the restriction $V|_M$ is nowhere zero. Define the index $\ind_x V$ of $V$ at an isolated singularity $x$ as usual (see, e.g., \cite[p. 136]{hirsch}), and let $\ind V=\sum_{x\in \sing V} \ind_x V$ 
denote the sum of the local indices. 
Also define
\begin{equation}
\label{eq:alpha}
\alpha_V:M\to STX|_M;\ x\mapsto \frac{V(x)}{|V(x)|}.
\end{equation}
by rescaling $V$. 

Following~\cite{sha}, we assume throughout the paper the following condition:
\ben\label{lp}
\text{the metric on }X\text{ is \emph{locally product} near the boundary }M,
\een
which in particular implies that $M$ is a totally goedesic submanifold of $X$. The general case is addressed in \cite{nie}. 

Under condition \er{lp}, Sha~\cite{sha} proved his version of the relative Poincar\'e-Hopf theorem
\begin{equation}
\label{sha's ph}
\ind V-\int_{\alpha_V(M)}\Phi=\begin{cases}
\chi(X) & \text{ if }\dim  X\text{ is even,}\\
0 &  \text{ if }\dim  X\text{ is odd.}
\end{cases}
\end{equation}

The starting point of this paper is to study $\Phi$, or rather its certain restriction defined as follows. Let $\vn$ denote the outward unit normal vector field of $M$. The images $\vec n(M)$ and $(-\vec n)(M)$ in $STX|_M$ 
are the spaces of outward and inward unit normal vectors of $M$. Define
\begin{equation}
\label{eq-def-cstm}
CSTM:=STX|_M\backslash(\vec n(M)\cup (-\vec n)(M))
\end{equation}
 ($C$ for cylinder) to be the complement.

\begin{theorem}
\label{th-exact}
Under condition \er{lp}, $\Phi$ is exact on $CSTM$ \er{eq-def-cstm}. More precisely, there is a differential form $\Gamma$ of degree $n-2$ on $CSTM$ such that
\begin{equation*}
\Phi=d\Gamma.
\end{equation*}
\end{theorem}

The definition of $\Gamma$ is in Definition \ref{def-gamma}, and the above theorem is proved right after that. 

Theorem \ref{th-exact} and Stokes' theorem then allow the following concrete evaluation of Sha's term $\int_{\alpha_V(M)}\Phi$ in \er{sha's ph} in terms of more classical local indices. For a generic vector field $V$, let $\partial V$ be the projection of $V|_M$ to $TM$ according to \er{collar}, and let 
$\partial_- V$ (resp. $\partial_+ V$) be the restriction of $\partial V$ to the subspace of $M$ where $V$ points inward (resp. outward) to $X$. Generically $\partial_\pm V$ have isolated singularities. (A non-generic $V$ can always be modified by adding an extension to $X$ of a normal vector field or a tangent vector field to $M$.)

\begin{theorem} 
\label{th-index}
Under condition \er{lp} and for a generic vector field $V$, one has
\begin{equation}
\label{eq-index}
\int_{\alpha_V(M)}\Phi=
\begin{cases}
-\ind\pmv & \text{ if }\dim X \text{ is even,}\\
\frac{1}{2} (\ind \ppv-\ind\pmv) & \text{ if }\dim X \text{ is odd}.
\end{cases}
\end{equation}
\end{theorem}

\begin{remark}\label{misc}
Generically $\ind\ppv+\ind\pmv=\ind\partial V=\chi(M)$ by the Poincar\'e-Hopf theorem. When $\dim X$ is even and hence $\dim M$ is odd, since $\chi(M)=0$, one has equality between the two formulas in \er{eq-index}. 

When $\dim X$ is odd, since $\chi(M)=2\chi(X)$ by basic topological knowledge, one has the following reformulation of the odd case in \er{eq-index}
\ben\label{reform}
\int_{\alpha_V(M)}\Phi=\frac 1 2 \chi(M)-\ind\pmv=
\chi(X)-\ind\pmv,\ \text{if }\dim X\text{ is odd}.
\een
\end{remark}

We finish this introduction by explaining the relation of our result with the Law of Vector Fields. For a generic vector field $V$, using 
the flow along $-V$ and counting fixed points with multiplicities, one has the following 
\emph{Law of Vector Fields}:
\ben
\label{law of v f}
\ind V+\ind \pmv=\chi(X).
\een
This was first proved by Morse \cite{morse} and later on publicized by Gottlieb, who also coined the term. 

Our result \er{eq-index} and the reformulation \er{reform}  of the odd case 
then  directly show that the two relative Poincar\'e-Hopf theorems, \er{sha's ph} and \er{law of v f}, are equivalent. Therefore following the route of the relative Poincar\'e-Hopf theorem of Sha \cite{sha} under condition \er{lp}, our result \er{eq-index} gives a purely differential-geometric proof of the Law of Vector Fields. Other differential-geometric proofs are given in \cite{nie}. 

\smallskip

\noindent\emph{Acknowledgement.} The author would like to thank Wojciech Dorabia\l a for getting him interested in  this topic, in particular Sha's paper \cite{sha}, and for motivating discussions. He also thanks the referee for careful reading and comments. 

\section{Differential forms}

Throughout the paper, $c_{j-1}$ denotes the volume of the unit $(j-1)$-sphere $S^{j-1}$. 

Chern's transgression form $\Phi$ is defined as follows. 
Choose oriented local orthonormal frames $\{e_1,e_2,\cdots,e_n\}$ for the tangent bundle $TX$. Let $(\omega_{ij})$ and 
$(\Omega_{ij})$ be the $\mathfrak{so}(n)$-valued connection forms and curvature forms
for the Levi-Civita connection $\nabla$ of the Riemannian metric on $X$ defined by 
\begin{gather}
\nabla e_i=\sum_{k=1}^n \omega_{ij}e_j,\label{def-omega}\\
\Omega_{ij}=d\omega_{ij}-\sum_{k=1}^n \omega_{ik}\omega_{kj}.\label{def-O}
\end{gather}
Let the $u_i$ be the coordinate functions  on $STX$ in terms of the frames defined by
\ben
\label{def-u}
v=\sum_{i=1}^{n} u_i(v)e_i,\quad \forall v\in STX.
\een
Let the $\theta_i$ be the 1-forms on $STX$ defined by
\ben
\label{def-theta}
\theta_i=du_i+\sum_{k=1}^{n}u_k\omega_{ki}.
\een
For $ k=0,1,\cdots,[\frac{n-1}2]$ (with $[-]$ standing for the integral part), define the degree $n-1$ forms on $STX$
\ben
\label{eq-phi-j}
\Phi_k=\sum_{\tau} \epsilon(\tau)u_{\tau_1}\theta_{\tau_2}\cdots\theta_{\tau_{n-2k}}\Omega_{\tau_{n-2k+1}\tau_{n-2k+2}}\cdots \Omega_{\tau_{n-1}\tau_{n}},
\een
where the summation runs over all permutations $\tau$ of $\{1,2,\cdots,n\}$, and $\epsilon(\tau)$ is the sign of $\tau$. (The index $k$ stands for the number of curvature forms involved. Hence the restriction $0\leq k\leq [\frac{n-1}2]$. This convention applies throughout the paper.) Define Chern's transgression form as
\begin{align}
\Phi&=\frac{1}{(n-2)!!c_{n-1}} \sum_{k=0}^{[\frac{n-1}2]} (-1)^k \frac{1}{2^k k! (n-2k-1)!!} \Phi_k\label{phi}\\
&=:\frac{1}{(n-2)!!c_{n-1}} \sum_{k=0}^{[\frac{n-1}2]} \overline{\Phi_k}=:\sum_{k=0}^{[\frac{n-1}2]} \wt{\Phi_k}.\nm
\end{align}
(See \er{eq-basic} for an explanation, in the case of $M$ with dimension $n-1$, for the coefficients involved.)
The $\Phi_k$ and hence $\Phi$ are invariant under $\mathrm{SO}(n)$-transformations of the local frames and hence are intrinsically defined. Note that the 0th term  
\begin{equation}
\label{0th term}
\wt{\Phi_0}=\frac{1}{(n-2)!!c_{n-1}} \frac 1{(n-1)!!}\Phi_0=\frac{1}{c_{n-1}}d\sigma_{n-1}=\wt{\vol}_{n-1}
\end{equation}
is the relative unit volume form of the fibration $S^{n-1}\to STX\to X$, since by \er{eq-phi-j}
\ben\label{factorial}
\Phi_0=\sum_\tau \epsilon(\tau)u_{\tau_1}\theta_{\tau_2}\cdots\theta_{\tau_n}=(n-1)!d\sigma_{n-1}
\een
(see \cite[(26)]{chern}). 

Now we start to transgress $\Phi$ \er{phi} on $CSTM$ \er{eq-def-cstm}. At $TX|_M$, we choose oriented local orthonormal frames $\{e_1,e_2,\cdots,e_n\}$ such that $e_1=\vec n$ is the outward unit normal vector of $M$. Therefore $\{e_2,\cdots,e_{n}\}$ are oriented local orthonormal frames for $TM$. 
Let $\phi$ be the angle coordinate on $STX|_M$ defined by
\ben
\label{def-phi}
\phi(v)=\angle(v,e_1)=\angle(v,\vn),\ \forall v\in STX|_M. 
\een
One has from \er{def-u}
\begin{equation}
\label{eq-u1}
u_1=\cos\phi.
\end{equation}
Let
\begin{align}
p:CSTM=STX|_M\backslash(\vn(M)\cup (-\vn)(M))\to STM;\ v&\mapsto \frac{\partial v}{|\partial v|}\nonumber\\
\text{(in coordinates)} (\cos\phi,u_2,\cdots,u_{n})&\mapsto\frac{1}{\sin\phi}(u_2,\cdots,u_{n})
\label{def p}
\end{align}
be the projection to the equator $STM$.
By definition, 
\begin{gather}
\text{for }x\in M,\ \partial V(x)=0\Leftrightarrow \a_V(x)=\pm \vec n(x), 
\label{equal}\\
\label{pa}
p\circ \a_V=\a_{\partial V}\text{ when }\partial V\neq 0. 
\end{gather}

The locally product metric \er{lp} near $M$ means that $\nabla e_1=\nabla\vec n=0$. Hence from \er{def-omega} one has
\ben
\label{o=0}
\omega_{1*}=-\omega_{*1}=0.
\een 
From \er{def-theta}, \er{eq-u1} and \er{o=0}, one has
\ben
\label{eq-theta1}
\theta_1=-\sin\phi\,d\phi.
\een 
From \er{def-O} and \er{o=0}, one also has 
\begin{equation}
\label{eq-1*=0}
\Omega_{1*}=-\Omega_{*1}=0.
\end{equation}

We use the convention that $\tau$ is a permutation of $(1,2,\cdots,n)$ and $\rho$ is a permutation of $(2,\cdots,n)$. 

In view of \eqref{eq-1*=0} on $STX|_M$, the index 1 in the formula \eqref{eq-phi-j} for $\Phi_k$ appears in either $u_{\tau_1}$ or one of the $\theta_{\tau_i}$ for $2\leq i\leq n-2k$. There are totally $n-2k-1$ possibilities for the second case. 

Therefore, on $STX|_M$, one has the following more concrete 
\begin{equation}
\label{eq-decomp}
\Phi_k=u_1 \Xi_k-(n-2k-1)\theta_1 \Upsilon_k,\ k=0,\cdots,\m,
\end{equation}
where 
\begin{gather}
\Upsilon_k=\sum_{\rho} \epsilon(\rho)u_{\rho_2}\theta_{\rho_3}\cdots\theta_{\rho_{n-2k}}\Omega_{\rho_{n-2k+1}\rho_{n-2k+2}}\cdots \Omega_{\rho_{n-1}\rho_{n}}, \ k=0,\cdots,[\frac {n-2}2],
\label{eq-phi-i}\\
\Xi_k=\sum_{\rho} \epsilon(\rho)\theta_{\rho_2}\theta_{\rho_3}\cdots\theta_{\rho_{n-2k}}\Omega_{\rho_{n-2k+1}\rho_{n-2k+2}}\cdots \Omega_{\rho_{n-1}\rho_{n}},\ k=0,\cdots,[\frac {n-1}2].
\label{eq-lambda}
\end{gather}
The negative sign in \er{eq-decomp} is from $\epsilon(\tau)$ in \er{eq-phi-j} when one moves $\theta_1$ in front of $u_{\tau_1}$. 
(When $[\frac {n-1}2]=[\frac {n-2}2]+1$, one can either define $\Upsilon_{[\frac {n-1}2]}=0$ by dimensional reason or observe its coefficient in \er{eq-decomp} $n-2k-1=0$ for $k=[\frac {n-1}2]$. This observation applies throughout the section.)

We use the convention to write superscript ${}^e$ for functions and forms defined on the equator $STM$ of $STX|_M$. 
Using the $\{e_2,\cdots,e_{n}\}$ as oriented local orthonormal frames for $TM$, we define $u^e_i$, $\theta^e_i$ and $\Phi^e_k$ as functions and forms on $STM$ in the same way as 
in \er{def-u}, \er{def-theta} and \er{eq-phi-j}. Note that $\Omega^e_{ij}=\Omega_{ij}$ for $2\leq i,j\leq n$ by \er{def-O} and \er{o=0}. Therefore one has the degree $n-2$ forms on $STM$
\begin{gather}
\Phi^e_k=\sum_{\rho} \epsilon(\rho)u^e_{\rho_2}\theta^e_{\rho_3}\cdots\theta^e_{\rho_{n-2k}}\Omega_{\rho_{n-2k+1}\rho_{n-2k+2}}\cdots \Omega_{\rho_{n-1}\rho_{n}},\ k=0,\cdots,[\frac {n-2}2].
\label{eq-phi-i-e}
\end{gather}
Following \cite{chern}, also define the degree $n-1$ forms on $STM$
\ben
\Psi^e_k=\sum_{\rho} \epsilon(\rho)\theta^e_{\rho_2}\theta^e_{\rho_3}\cdots\theta^e_{\rho_{n-2k}}\Omega_{\rho_{n-2k+1}\rho_{n-2k+2}}\cdots \Omega_{\rho_{n-1}\rho_{n}},\ k=0,\cdots,[\frac {n-1} 2].
\label{eq-lambda-e}
\een
Note that the $\Phi_k^e$ and the $\Psi_k^e$ are just the $\Upsilon_k$ in \er{eq-phi-i} and the $\Xi_k$ in 
\er{eq-lambda} with the superscript ${}^e$. 
By dimensional  reason one has 
\begin{gather}
\Psi^e_0=0.
\label{eq-lamb0e}
\end{gather} 
On $STM$, Chern's basic formulas~\cite{chern} are 
\begin{equation}
\label{eq-basic}
d\Phi^e_k=\Psi^e_k+\frac{n-2k-2}{2(k+1)}\Psi^e_{k+1},\ k=0,\cdots,[\frac {n-2}2].
\end{equation}
(This also explains, over $M$ with dimension $n-1$, the construction of $\Phi$ in \er{phi} for the purpose of consecutive cancellations.) 

\begin{lemma} One has on $CSTM$ \er{eq-def-cstm}, for $k=0,\cdots,[\frac {n-1} 2]$, 
\begin{equation}
\label{eq-psij}
\Phi_k=\sin^{n-2k-1}\phi\cos\phi\, p^*\Psi^e_k+(n-2k-1)\sin^{n-2k-2}\phi\,d\phi\, p^*\Phi^e_k.
\end{equation}
\end{lemma}

\begin{proof}
For $2\leq i\leq n$ and from \er{def p}, one has
\begin{equation}
p^*u^e_i=\frac{1}{\sin\phi} u_i.
\label{eq-p&u}
\end{equation}
Differentiating the above and using \er{def-theta} and \er{o=0}, one has 
\begin{equation}
p^*\theta^e_i=\frac{1}{\sin\phi} \theta_i-\frac{\cos\phi}{\sin^2\phi}d\phi\, u_i
\label{eq-p&theta}
\end{equation}
Because of the presence of $d\phi$ and
in view of \er{eq-phi-i-e}, \er{eq-p&u} and \er{eq-p&theta}, one has
\begin{gather}
\label{eq-t1phij}
d\phi\, p^*\Phi^e_k=\frac{1}{\sin^{n-2k-1}\phi}\,d\phi\,\Upsilon_k\Rightarrow \,d\phi\,\Upsilon_k=\sin^{n-2k-1}\phi\, d\phi\, p^*\Phi^e_k,
\end{gather}
where $n-2k-1$ is the number of $u$ and $\theta$'s in \er{eq-phi-i-e}. 
Hence by \er{eq-theta1}, one has
\begin{gather}
\label{using theta1}
\theta_1 \Upsilon_k=-\sin^{n-2k}\phi\,d\phi\,p^*\Phi^e_k,
\end{gather}

Now the pullback of $\Psi^e_k$ in \er{eq-lambda-e} is slightly harder, since $d\phi$ may come up as in \er{eq-p&theta}, but only once among the $(n-2k-1)$ $\theta^e$'s. Therefore
\begin{equation}
p^*\Psi^e_k=\frac{1}{\sin^{n-2k-1}\phi} \Xi_k-(n-2k-1) \frac{\cos\phi}{\sin^{n-2k}\phi} d\phi\, \Upsilon_k.
\end{equation}
Using \er{eq-u1} and \er{eq-t1phij}, one then has
\begin{align}
u_1\Xi_k&=\sin^{n-2k-1}\phi\cos\phi\, p^*\Psi^e_k+(n-2k-1) \frac{\cos^2\phi}{\sin\phi}d\phi\, \Upsilon_k\label{eq-u1lambj}\\
&=\sin^{n-2k-1}\phi\cos\phi\, p^*\Psi^e_k+(n-2k-1) \sin^{n-2k-2}\phi\cos^2\phi \,d\phi\, p^*\Phi^e_k.\nm
\end{align}

Combining \er{eq-decomp}, \er{using theta1} and \er{eq-u1lambj}, one has
\begin{align*}
\Phi_k=&u_1\Xi_k-(n-2k-1)\theta_1 \Upsilon_k\\
=&\sin^{n-2k-1}\phi\cos\phi\, p^*\Psi^e_k+(n-2k-1) \sin^{n-2k-2}\phi\cos^2\phi \,d\phi\, p^*\Phi^e_k\\
&+(n-2k-1)\sin^{n-2k}\phi\, d\phi\, p^*\Phi^e_k\\
=&\text{RHS of }\er{eq-psij}
\end{align*}
by $\cos^2\phi+\sin^2\phi=1$. 
\end{proof}

Since $\Psi^e_0=0$ \er{eq-lamb0e}, one has from \er{eq-psij}
\begin{equation}
\label{eq-psi0}
\Phi_0=(n-1)\sin^{n-2}\phi\,d\phi\, p^*\Phi^e_0.
\end{equation}

\begin{remark}
In view of \er{factorial}, \er{eq-psi0} is just the relation (due to condition \er{lp}) between the relative volume forms $d\sigma_{n-1}$ of $S^{n-1}\to STX|_M\to M$ and $d\sigma_{n-2}$ of $S^{n-2}\to STM\to M$
$$
d\sigma_{n-1}=\sin^{n-2}\phi\,d\phi\,p^*d\sigma_{n-2}.
$$
On one fixed sphere and its equator, this is an easy fact and follows from using spherical coordinates, which also
accounts for the basic formula
\ben\label{basic-know}
c_{n-1}= {c_{n-2}}\int_{0}^\pi \sin^{n-2}\phi\,d\phi.
\een 
\end{remark}

Our goal is to find a differential form $\Gamma$ such that $d\Gamma=\Phi$. We do this inductively starting from the above $\Phi_0$ in \er{eq-psi0}. Therefore we need to use an antiderivative of $\sin^{n-2}\phi$. 

\begin{definition}
\label{def ik}
For a non-negative integer $b$, define functions of $\phi$
\begin{equation}
\label{eq-ik}
I_b(\phi)=\int \sin^b \phi\,d\phi,
\end{equation}
where we require the arbitrary constants to be 0. More precisely, 
\begin{equation}
\label{different}
I_b(\phi)=\begin{cases}
\int_0^\phi \sin^b t\,dt, & \text{ if }b\text{ is even},\\
\int_{\frac{\pi}{2}}^\phi \sin^b t\,dt, & \text{ if }b\text{ is odd}.
\end{cases}
\end{equation}

Integration by parts gives
\begin{equation}
\label{eq-ikind}
bI_b(\phi)+\sin^{b-1}\phi\cos\phi=(b-1)I_{b-2}(\phi),\ b\geq 2.
\end{equation}
Clearly $I_0(\phi)=\phi$ and $I_1(\phi)=-\cos\phi$. These also inductively determine $I_b(\phi)$. 
\end{definition}

\begin{definition}\label{def-gamma}
 We define the following differential forms of degree $n-2$ on $CSTM$ \er{eq-def-cstm} 
\begin{gather}
\label{eq-Gamj}
\Gamma_k=I_{n-2k-2}(\phi) p^*\Phi^e_k,\ k=0,\cdots,[\frac {n-2} 2],\\
\label{eq-tgamm}
\overline{\Gamma_k}=(-1)^k\frac{1}{2^k k! (n-2k-3)!!} \Gamma_k=(-1)^k\frac{1}{2^k k! (n-2k-1)!!} (n-2k-1)\Gamma_k
\end{gather}
(with the convention $(-1)!!=1$), and 
\begin{align}
\label{gamma}
\Gamma&=\frac{1}{(n-2)!!c_{n-1}}\sum_{k=0}^{[\frac {n-2} 2]} \overline{\Gamma_k}\\
&=\frac{1}{(n-2)!!c_{n-1}}\sum_{k=0}^{[\frac {n-2} 2]}(-1)^k\frac{1}{2^k k! (n-2k-3)!!} I_{n-2k-2}(\phi) p^*\Phi^e_k .\nm
\end{align}
\end{definition}

With this definition of $\Gamma$, now we prove Theorem \ref{th-exact}. 

\begin{proof}[Proof of Theorem \ref{th-exact}]
First by Chern's basic formula \er{eq-basic}, for $k=0,\cdots,[\frac {n-2} 2]$,
\ben
\label{dgam}
d\Gamma_k=\sin^{n-2k-2}\phi \,d\phi\, p^*\Phi^e_k+I_{n-2k-2}(\phi)p^*\Psi^e_k+\frac{n-2k-2}{2(k+1)} I_{n-2k-2}(\phi)p^* \Psi_{k+1}^e.
\een

Define 
\ben
\label{def lj}
\overline{L_k}= (-1)^k\frac{n-2k}{2^k k! (n-2k-1)!!} I_{n-2k}(\phi)p^*\Psi^e_k,\ k=0,\cdots,[\frac {n} 2]
\een
($L$ for leftover). 
 
\begin{claim} For $k=0,\cdots,[\frac {n-2} 2]$, one has
\ben
\label{p-dg=l}
\sum_{i=0}^k \overline{\Phi_i}-d(\sum_{i=0}^k\overline{\Gamma_i})=\overline{L_{k+1}}.
\een
\end{claim}

\begin{proof}[Proof of the claim]
By induction. Actually \er{p-dg=l} clearly holds for $k=-1$, since both sides are zero by natural reasons, \er{def lj} and \er{eq-lamb0e}. 
(One can also check the $k=0$ case using the same reason as in the following induction step.) 

Now assume \er{p-dg=l} holds for $k-1$. Then using this induction hypothesis, plugging in all the formulas \er{def lj}, \er{dgam}, \er{eq-tgamm}, \er{eq-psij}, and by \er{eq-ikind}, one has
\begin{align*}
&\sum_{i=0}^k \overline{\Phi_i}-d(\sum_{i=0}^k\overline{\Gamma_i})=\overline{L_k}+\overline{\Phi_k}-d\overline{\Gamma_k}\\
=&(-1)^k\frac{1}{2^k k! (n-2k-1)!!} [((n-2k)I_{n-2k}(\phi)+\sin^{n-2k-1}\phi\cos\phi)p^*\Psi^e_k\\
&+(n-2k-1)\sin^{n-2k-2}\phi\,d\phi\,p^*\Phi_k^e-(n-2k-1)\sin^{n-2k-2}\phi \,d\phi\, p^*\Phi^e_k\\
&-(n-2k-1)I_{n-2k-2}(\phi)p^*\Psi^e_k-\frac{(n-2k-2)}{2(k+1)}(n-2k-1) I_{n-2k-2}(\phi)p^* \Psi_{k+1}^e]\\
=&(-1)^{k+1}\frac{n-2k-2}{2^{k+1} (k+1)! (n-2k-3)!!}I_{n-2k-2}(\phi)p^*\Psi_{k+1}^e=\overline{L_{k+1}}.
\end{align*}
\end{proof}

When $n=2m$ for $m\geq 1$, $[\frac {n-2} 2]=[\frac {n-1} 2]=m-1$. Therefore to prove $\Phi=d\Gamma$, in view of \er{phi} and \er{gamma}, it suffices by \er{p-dg=l} to proceed as follows 
$$
\sum_{i=0}^{m-1} \overline{\Phi_i}-d(\sum_{i=0}^{m-1}\overline{\Gamma_i})=\overline{L_{m}}=0
$$
since $\overline{L_m}=0$ from \er{def lj} due to the coefficient $n-2k$ on the top. 

When $n=2m+1$ for $m\geq 1$, $[\frac {n-2} 2]=m-1$ and $[\frac {n-1} 2]=m$.  
In view of \er{phi}, \er{gamma}, \er{p-dg=l}, \er{def lj} and \er{eq-psij}, one has
\begin{align}
 &\sum_{i=0}^{m} \overline{\Phi_i}-d(\sum_{i=0}^{m-1}\overline{\Gamma_i})=\overline{L_{m}}+\overline{\Phi_m}\nm\\
 =&(-1)^m \frac{1}{2^m m!}(I_1(\phi)+\cos\phi)p^*\Psi_m^e=0,\label{=0?}
 \end{align}
 since $I_1(\phi)=-\cos\phi$ by Definition \ref{def ik}. The proof is now complete.
 \end{proof}

\section{Indices}

Now we are ready for the proof of Theorem \ref{th-index} using Stokes' theorem. 

\begin{proof}[Proof of Theorem \ref{th-index}] 
Let $B_r^M(\sing \partial V)$ (resp. $S_r^M(\sing \partial V))$ denote the union of smalls open balls (resp. spheres) of radii $r$ in $M$ around the finite set of points $\sing \partial V$. 
Then by \er{equal}, $\a_V(M\backslash B^M_r(\sing \partial V))\subset CSTM$. 
By Theorem \ref{th-exact} and Stokes' theorem, 
\begin{align}
\label{stokes 1}
&\int_{\alpha_V(M)}\Phi=\lim_{r\to 0} \int_{\a_V(M\backslash B_r^M(\sing\partial V))}\Phi
= \lim_{r\to 0} \int_{\a_V(M\backslash B_r^M(\sing\partial V))} d\Gamma\\ \nonumber
=&-\lim_{r\to 0} \int_{\a_V(S_r^M(\sing \partial V))} \Gamma
=-\lim_{r\to 0} \int_{\a_V(S_r^M(\sing\partial V))}  \frac{1}{(n-2)!!c_{n-1}}\overline{\Gamma_0},
\end{align}
since all the other $\overline{\Gamma_k}$ for $k\geq 1$ in \er{gamma} involve curvature forms and hence don't contribute in the limit when integrated over small spheres (see \cite[\S 2]{chern2}). 

One has by Definition \ref{def-gamma}
\begin{equation}
\label{vol again}
\frac{1}{(n-2)!!c_{n-1}}\overline{\Gamma_0}=\frac{1}{(n-2)!!c_{n-1}}\frac{1}{(n-3)!!} I_{n-2}(\phi)p^*\Phi_0^e=\frac 1 {c_{n-1}} I_{n-2}(\phi)p^*d\sigma_{n-2}
\end{equation}
with $\vol_{n-2}$ being the relative volume form of $S^{n-2}\to STM\to M$, since $\Phi_0^e=(n-2)!\vol_{n-2}$ (cf. \er{factorial}).

Continuing \er{stokes 1} and using \er{vol again}, one has
\begin{align*}
& \int_{\alpha_V(M)}\Phi=-\frac{1}{c_{n-1}} \lim_{r\to 0} \int_{\a_V(S_r^M(\sing\ppv)\cup S_r^M(\sing \pmv))} I_{n-2}(\phi)p^*d\sigma_{n-2}\\
\overset{(1)}=& -\frac{1}{c_{n-1}} [I_{n-2}(0)\lim_{r\to 0} \int_{\a_{\partial V}(S_r^M(\sing\ppv))} \vol_{n-2}\\
&+ I_{n-2}(\pi)\lim_{r\to 0} \int_{\a_{\partial V}(S_r^M(\sing\pmv))} \vol_{n-2}]\\
\overset{(2)}=&-\frac{c_{n-2}}{c_{n-1}} (I_{n-2}(0)\ind\ppv+I_{n-2}(\pi)\ind\pmv)\\
\overset{(3)}=&\begin{cases}
-\ind\pmv & \text{ if }n=\dim X\text{ is even},\\
\frac{1}{2} (\ind\ppv-\ind\pmv) & \text{ if }n=\dim X\text{ is odd}.
\end{cases}
\end{align*}
Here equality (1) uses \er{pa} and
\begin{gather*}
\phi(\alpha_V(x))\to \pi \text{ for }x\in S_r^M(\sing\pmv),\text{ as }r\to 0,\\
\phi(\alpha_V(x))\to 0 \text{ for }x\in S_r^M(\sing\ppv),\text{ as }r\to 0.
\end{gather*}
Equality (2) is by the definition of index. In view of \er{different}, one has
\begin{gather}
I_{n-2}(0)=0,\ I_{n-2}(\pi)=\int_0^\pi \sin^{n-2}\phi\,d\phi,\text{ if }n\text{ is even},\label{even for both}\\
I_{n-2}(0)=-\frac 1 2\int_0^\pi \sin^{n-2}\phi\,d\phi,\ I_{n-2}(\pi)=\frac 1 2\int_0^\pi \sin^{n-2}\phi\,d\phi,\text{ if }n\text{ is odd},\label{odd symmetry}
\end{gather}
where \er{odd symmetry} uses symmetry of integrals. 
Then equality (3) follows from \er{basic-know}.
\end{proof}

\begin{remark} If instead of \er{different}, one also defines ${I_b}(\phi)=\int_0^\phi \sin^b t\, dt$ for the odd case, it can be checked that one gets formulas different from but equivalent to ours. 
\end{remark}

\end{document}